\documentclass[onecolumn]{autart}   

\usepackage{amsmath,amssymb} 
\usepackage{amssymb}  
\usepackage{tikz}
\usepackage{multicol}
\usepackage{blindtext}
\newtheorem{theorem}{\textbf{Theorem}}

\def\E{\mathbf{E}}

\begin{document}

\begin{frontmatter}

\title{On Optimal Distributed Output-Feedback Control over Acyclic Graphs}

\author{
Ather Gattami
},
\author{
Omid Khorsand
}

\address{
Electrical Engineering School, KTH-Royal Institute of Technology, SE-100 44 Stockholm, Sweden.
}
\end{frontmatter}

\begin{abstract}
In this paper, we consider the problem of distributed optimal control of linear dynamical systems 
with a quadratic cost criterion. We study the case of output feedback control for 
two interconnected dynamical systems, and show that the linear optimal solution can be obtained from
a combination of two uncoupled Riccati equations and two coupled Riccati equations. 
\end{abstract}
\section{INTRODUCTION}
Decentralized control has been a challenging area of open research for over three decades. However, these problems are still intractable in general. In one of the earliest works, it was shown that for a linear system subject to quadratic cost and Gaussian noise, complex nonlinear controllers could outperform any given linear controller (see \cite{c1}).

Centralized controller design problems assume a classical information pattern where each system has access to the information of all other systems. However, in decentralized problems, different systems have access to different information sets and their decisions are based on local measurements. Solving the problem in the presence of this constraint is a much harder problem.

One approach has been to classify specific information patterns leading to linear optimal controllers. In \cite{c4}, sufficient conditions are given under which optimal controllers are linear for the LQG problem. The one-step delay LQG problem was solved in \cite{athans}. An important result was given in \cite{c5} which showed that for a new information structure that was called \emph{partially nested}, the optimal controller is linear in the information set. In \cite{c6}, stochastic linear quadratic control problem was solved under the condition that all the subsystems have access to the global information from some time in the past. In \cite{bamie}, it was shown that the constrained linear optimal decision problem for infinite horizon linear quadratic control, can be posed as an infinite dimensional convex optimization problem, given that the system considered is stable. A solution to the state feedback problem, using a duality result for distributed estimation and control under partially nested information pattern, was given in \cite{c81} .

In this paper, we consider a decentralized system over an acyclic graph. A quadratic cost is to be minimized under the partially nested information structure. The goal is to explore the structure of the distributed controllers.

This problem is known to have a linear optimal policy, see \cite{c5} and \cite{c9}. However, most existing approaches do not provide explicit optimal controller formulae and, in addition, the order of the controller can be large. A solution to this problem has been first given in \cite{thesis}, where the resulting controller is of high order which makes the implementation difficult. Moreover, some work has been focused on finding numerical algorithms to these problems \cite{num1}, \cite{num2}.

State-space solutions to the two-systems \emph{state feedback} $\mathcal{H}_2$ version of this problem was given in \cite{c12}. Another instance of the problem  with a \emph{partial output feedback} structure has been considered in \cite{c13}, in which one system measures his state directly and the other system uses a noisy measurement of his own state. However, no previous work has provided a solution to the complete \emph{output feedback} problem for $N$ interconnected systems, not even for the case $N=2$. In this paper, we propose a methodology to attack the problem for the case of two interconnected systems to understand the structure of the controller.


\subsection*{Notation}
We denote a matrix partitioned into blocks by $A=[A_{ij}]$, where $A_{ij}$ denotes the block matrix of $A$ in block position $(i,j).$ A discrete-time stochastic process $x(0), x(1),x(2),...$ is denoted by $\{x(t)\}$. For vectors $y(0), y(1),...,y(t)$, we let $y(0:t)=\begin{bmatrix}y^T(0) &y^T(1)&\cdots&y^T(t)\end{bmatrix}^T$. $\delta(t)$ denotes the Dirac delta function.
\section{Problem Statement} \label{sec2}
Consider a linear system composed of $N$ interconnected subsystems
\begin{align*} x(t+1)&=Ax(t)+Bu(t)+w(t), \\
              y(t)&=Cx(t)+v(t), \end{align*}
where the system matrices are partitioned into blocks as $A=[A_{ij}],~B=[B_{ij}],~C=[C_{ij}],~i,j=1,...,N$, and $x$, $u$, and $y$ are the overall state, input, and output vectors. These vectos are partitioned as
\begin{align*}  x(t)&=\begin{bmatrix}x^T_{1}(t)&x^T_2(t)& \cdots &x^T_N(t) \end{bmatrix}^T, \\
                u(t)&=\begin{bmatrix}u^T_{1}(t)&u^T_2(t)& \cdots &u^T_N(t) \end{bmatrix}^T,\\
                y(t)&=\begin{bmatrix}y^T_{1}(t)&y^T_2(t)& \cdots &y^T_N(t) \end{bmatrix}^T,
\end{align*}
where $x_i(t)$ is the state of subsystem $i$, $y_i(t)$ is the noisy measurement made by subsystem $i$ and $u_i(t)$ is the control input to the subsystem $i$ at time $t$.
The interconnection structure can be described as follows. If the state of system $j$ at time step $t$ (i.e., $x_j(t)$) affects the state of system $i$ at time step $t + 1$ (i.e., $x_i(t + 1)$), then $A_{ij} \neq 0$, otherwise $A_{ij} = 0$. This structure can be represented by a directed graph of order $N$. The graph has an arrow from node $j$ to $i$ if and only if $A_{ij}\neq 0$.

Here, we assume system dynamics over an acyclic graph. For a directed acyclic graph it is natural to order the nodes so that the adjacency matrix becomes lower-triangular. Thus, we shall consider a lower-triangular $A$ matrix for the remainder of this paper. Furthermore, we assume that the matrices  $B$ and $C$ have the same sparsity structure:
\begin{align}A&=\begin{bmatrix}A_{11} &    0    &    \cdots &  0      \\
                              A_{21} & A_{22}  &    \cdots &  0      \\
                              \vdots & \vdots  &    \ddots & \vdots  \\
                              A_{N1} & A_{N2}  &    \cdots & A_{NN}
                              \end{bmatrix},\label{A} \\B&=\begin{bmatrix}B_{11} &    0    &    \cdots &  0      \\
                              B_{21} & B_{22}  &    \cdots &  0      \\
                              \vdots & \vdots  &    \ddots & \vdots  \\
                              B_{N1} & B_{N2}  &    \cdots & B_{NN}
                              \end{bmatrix}, \label{B} \\
                              C&=\begin{bmatrix}C_{11} &    0    &    \cdots &  0      \\
                              C_{21} & C_{22}  &    \cdots &  0      \\
                              \vdots & \vdots  &    \ddots & \vdots  \\
                              C_{N1} & C_{N2}  &    \cdots & C_{NN}
                              \end{bmatrix}, \label{C}
                              \end{align}
The interconnection structure for a two-systems chain, and a three-systems graph has been shown in Fig. \ref{fig1} and Fig. \ref{fig2} respectively.
Assume $\{w(t)\}$ and $\{v(t)\}$ are sequences of Gaussian vectors with zero mean values and covariances
\begin{align*}
&\mathbf{E}\begin{bmatrix}w(k)\\v(k)\end{bmatrix}\begin{bmatrix}w(l)\\v(l)\end{bmatrix}^T=\delta(k-l)\begin{bmatrix}W&U\\U^T&V\end{bmatrix}.
\end{align*}
Furthermore, the covariance matrices are partitioned into appropriately sized blocks
\begin{align}W=[W_{ij}], V=[V_{ij}], \mbox{~for } i, j=1, \hdots, N. \end{align}
$V_j$ will denote the $j$th block column of $V$. We want to minimize the infinite-horizon quadratic cost
\begin{equation} \label{cost}
 J= \lim_{M \rightarrow \infty} \frac{1}{M}  \sum_{t=0}^{M-1} \E \begin{bmatrix} x(t) \\ u(t) \end{bmatrix}^T
 \begin{bmatrix}Q&S\\S^T & R \end{bmatrix}
 \begin{bmatrix} x(t) \\ u(t) \end{bmatrix},
\end{equation}
where $\begin{bmatrix}Q&S\\S^T & R \end{bmatrix} \succ 0.$
The weighting matrices are also partitioned as $Q=[Q_{ij}], R=[R_{i j}], S=[S_{ij}],~i,j=1,...,N$.
$S_j$ will denote the $j$th block column of $S$.
In this system, the dynamics of subsystem $j$ propagates to subsystem $i$ if $j<i$. Clearly, if all systems have access to the global measurements, the information structure would be classical, and the optimal linear controller could be obtained from the LQG control theory. However, we consider the case in which each system $i$ has access to its own measurements up to time $t-1$, i.e. $y_i(0:t-1)$. Furthermore, system $i$ has access to the information set that system $j$ has access to if the decision that system $j$ makes affects the information set of system $i$. This assumption implies that the information pattern is \emph{partially nested}. This information pattern is not classical and is one of a few non-classical information patterns for which the optimal policy is known to be unique and linear in the associated information set \cite{c5}.

Hence, we are interested in finding the controllers restricted to the following structure:
\begin{align}
u_1(t)&=\mu_1(y_1(0:t-1)),            \nonumber      \\
u_2(t)&=\mu_2(y_1(0:t-1),y_2(0:t-1)), \nonumber      \\
&~\vdots                              \label{pattern}\\
u_N(t)&=\mu_N(y_1(0:t-1),\hdots, y_N(0:t-1)), \nonumber
\end{align}
where $\mu_1,\hdots,\mu_N$ denote linear functions in all their variables.
The objective is to find the stabilizing controllers that minimize the cost criterion (\ref{cost}) while restricted to the information constraints in (\ref{pattern}).
\section{Solution for Two Interconnected Systems}
\begin{figure}
\centering
\begin{tikzpicture}
\draw [->] (1,1) node[anchor=east,circle,
draw]{1} -- (2,1) node[anchor=west,
circle,draw]{2};
\end{tikzpicture}
\caption{Two-Systems Chain}
\label{fig1}
\vspace{1mm}
\begin{tikzpicture}
\draw [->] (0,0) node[anchor=east,circle,draw]{1} -- (1,0) node[anchor=west,circle,draw]{2};
\draw [->] (1.65,0) -- (2.5,0) node[anchor=west,circle,draw]{3};
\draw [<-]  (2.7,0.3) arc (50:131:2.2cm);
\end{tikzpicture}
\caption{Three-Systems Graph}
\label{fig2}
\end{figure}
For the two-systems problem, (Fig. \ref{fig1}), control inputs are constrained on the form:
\begin{align*}
u_1(t)&=\mu_1(y_1(0:t-1)),  \\
u_2(t)&=\mu_2(y_1(0:t-1),y_2(0:t-1)).
\end{align*}
Equivalently, controllers can be described as the following linear dynamical systems:
\begin{align*}
  \eta_1(t+1)&=F_1\eta_1(t)+G_1y_1(t)  \\
  u_1(t)&=-H_{1}\eta_1(t) 
\end{align*}
  \begin{align*}
  \eta_2(t+1)&=F_2\eta_2(t)+G_2\begin{bmatrix} y_1(t) \\ y_2(t) \end{bmatrix} \\
  u_2(t)&=-H_{2}\eta_2(t) 
  \end{align*}

To simplify notation, we define the block rows of $C$ as
$$\left[\begin{array}{c} C_{1} \\ \hline
                   {C}_{2}
\end{array}\right]=\left[
\begin{array}{cc} C_{11}   & 0 \\ \hline
                    C_{21} & C_{22}
\end{array}\right], $$
and the block columns of $B$ as

$$
\left[
\begin{array}{c|c} B_{1}   & B_2  
\end{array}\right]
=\left[
\begin{array}{c|c} 	B_{11} & 0 \\ 
                    		B_{21} & B_{22}
\end{array}\right].$$

We shall now state the optimal controller synthesis for the two-system problem as the main result of this section.
\subsection{\textbf{Main Result}}
\begin{theorem}
Consider the optimal control problem
\begin{align}\nonumber  \min_{u_1, u_2}&\lim_{M \rightarrow \infty} \frac{1}{M} \sum_{t=0}^{M-1} \E \begin{bmatrix} x(t) \\ u(t) \end{bmatrix}^T\begin{bmatrix} Q & S\\ S^T&R \end{bmatrix}\begin{bmatrix} x(t) \\ u(t) \end{bmatrix}     \\
\mbox{subject to~}
&\begin{bmatrix}
     x_1(t+1) \\
     x_2(t+1)
     \end{bmatrix}
  =
     \begin{bmatrix}
     A_{11}  & 0      \\
     A_{21}  & A_{22}
     \end{bmatrix}                 \begin{bmatrix}
                                    x_1(t) \\
                                    x_2(t)
                                   \end{bmatrix}      \nonumber      \\
 & \hspace{18 mm} +\begin{bmatrix}
  B_{11} & 0 \\
  B_{21}   & B_{22}
  \end{bmatrix} \begin{bmatrix}
  u_1(t) \\
  u_2(t)
  \end{bmatrix} +  \begin{bmatrix}
                  w_1(t) \\
                  w_2(t)
                  \end{bmatrix}                        \nonumber \\
 &\begin{bmatrix}
     y_1(t) \\
     y_2(t)
     \end{bmatrix}=
     \begin{bmatrix}
     C_{11} & 0 \\
     C_{21}      & C_{22}
     \end{bmatrix}                 \begin{bmatrix}
                                    x_1(t) \\
                                    x_2(t)
                                   \end{bmatrix} +  \begin{bmatrix}
                  v_1(t) \\
                  v_2(t)
                  \end{bmatrix}                    \nonumber   \\
&u_1(t)=\mu_{1}(y_1(0:t-1))  \nonumber \\
&u_2(t)=\mu_2(y_1(0:t-1),y_2(0:t-1)) \nonumber
\end{align}
Suppose there exist stabilizing solutions $P, \Pi, P_1, \Pi_2$, to the Riccati equations
\begin{align}
P =& (A-KC)P(A-KC)^T \nonumber \\
    +& W - KU - U^TK^T + KVK^T                \label{riccati4} \\
\Pi =& (A-BL)^T \Pi (A-BL) \nonumber \\
      +& Q - L^T S^T - S L + L^T R L                                           \label{riccati3} \\
P_1  = & 	(A-K_1C_1-B_2 L_2)P_1(A-K_1C_1-B_2 L_2)^T \nonumber \\
	+  &	 \mathcal{W} -K_1\mathcal {U}^T - \mathcal {U} K_1^T + K_1\mathcal{V} K_1^T		\label{riccati6}\\
\Pi_2=& (A-K_1C_1-B_2 L_2)^T\Pi_2 (A-K_1C_1-B_2 L_2)  \nonumber \\
	+ &  Q - L_2^T S_2^T - S_2 L_2 + L_2^T R_{22}L_2 \label{riccati5} 
\end{align}
where
 \begin{align}
           K&=(AP C^T+U)({C}P{C}^T+V)^{-1}, \label{K1}   \\
         L&=(B^T\Pi B+R)^{-1}(A^T\Pi B+S)^T, \label{L1}               \\
         K_1&=((A-B_2L_2)P_1{C}_1^T + \mathcal{U})({C}_1P_1{C}_1^T+\mathcal{V})^{-1}, \label{K2} \\
          L_2&=({B}_2^T\Pi_2{B}_2+R_{22})^{-1}((A-K_1C_1)^T\Pi_2{B}_2+S_2)^T, \label{L2} \\
         S_2&=\begin{bmatrix} S_{12} \\ S_{22} \end{bmatrix},                                                  \nonumber  
\end{align}
\begin{align}         
\begin{bmatrix} \mathcal{W}  & \mathcal{U} \\ \mathcal{U}^T& \mathcal{V}\end{bmatrix} &= \begin{bmatrix} K(CPC^T+V)K^T & KCPC_1^T + V_1 \\ C_1PC^TK^T+V_1^T & C_1PC_1^T+V_{11} \end{bmatrix}. 
\end{align}
Then, the {optimal controllers} are given by:
\begin{align} \label{ctrls1}
\begin{bmatrix} {u}_1(t) \\ u_2(t) \end{bmatrix} &= \begin{bmatrix} \hat{u}_1(t) \\ \hat{u}_2(t) \end{bmatrix} + 
\begin{bmatrix} 0 \\ \tilde{u}_2(t) \end{bmatrix}\\
\hat{z}(t+1) 	&= (A-K_1C_1)\hat{z}(t)+B\hat{u}(t)+K_1 y_1(t) \nonumber \\ 
\begin{bmatrix} \hat{u}_1(t) \\ \hat{u}_2(t) \end{bmatrix} 	&=  -L\hat{z}(t) \nonumber \\
z(t+1)  &= (A-KC) z(t) + Bu(t) + K y(t) \nonumber \\
\tilde{z}(t) 	   &=  {z}(t) - \hat{z}(t)	\nonumber \\
\tilde{u}_2(t)  &= -L_2 \tilde{z}(t). \nonumber 
 \end{align}		 

  \end{theorem}

\subsection{\textbf{Optimal Controller Derivation}}

\subsubsection{State and Controller Estimation}
In this part, we will present the optimal estimators involved in the controller design. Note that since the controllers $u_1$ and $u_2$ don't have the same information, their estimators of the system will be different.\\

Let 
\begin{align*}
\hat{x}(t) &= \E\{x(t) | y(0:t-1)\},\\
\tilde{x}(t) &= x(t) - \hat{x}(t).
\end{align*}
$\hat{x}(t)$ is given by the Kalman filter (\cite{c14})
\begin{equation*}
	\begin{aligned}
		\hat{x}(t+1)  &= A \hat{x}(t) + Bu(t) + K(C\tilde{x}(t)+v(t))
	\end{aligned}
\end{equation*}
where $K$ is obtained from (\ref{riccati4}) and (\ref{K1}). 
Introduce 
\begin{align*}
\omega(t) &= K(C\tilde{x}(t)+v(t)),\\
\nu(t)         &= C_1\tilde{x}(t)+v_1(t),\\
z(t) 		&= \hat{x}(t).
\end{align*}
By construction, $\nu(t)$ is orthogonal to $z(t)$.
Now we have that

\begin{align*}
\mathbf{E}\begin{bmatrix}\omega(t)\\ \nu(t)\end{bmatrix}\begin{bmatrix}\omega(t)\\ \nu(t)\end{bmatrix}^T &=\begin{bmatrix} \mathcal{W} & \mathcal{U} \\ \mathcal{U}^T& \mathcal{V}\end{bmatrix}\\
&=\begin{bmatrix} K(CPC^T+V)K^T & KCPC_1^T + V_1 \\ C_1PC^TK^T+V_1^T & C_1PC_1^T+V_{11} \end{bmatrix}
\end{align*}
and
\begin{equation*}
	\begin{aligned}
		z(t+1)  &= A z(t) + Bu(t) + \omega(t),\\
		y_1(t)  &= C_1z(t) + \nu(t).
	\end{aligned}
\end{equation*}
Also, let 
\begin{align*}
\hat{u}(t) &= \E\{u(t) | y_1(0:t-1)\},\\
\tilde{u}(t) &= u(t) - \hat{u}(t),\\
\hat{z}(t) &= \E\{z(t) | y_1(0:t-1)\}, \\
\tilde{z}(t) &= z(t) - z(t) -\hat{z}(t), 
\end{align*}
where $\hat{z}(t)$ is given by a Kalman filter
\begin{equation*}
	\begin{aligned}
		\hat{z}(t+1)  &= A \hat{z}(t) + B\hat{u}(t) + K_1(C_1\tilde{z}(t)+\nu(t)).\\
	\end{aligned}
\end{equation*}
and the estimation error dynamics are given by
\begin{align}\label{kalmanz2}
    \tilde{z}(t+1)=(A-K_1C_1)\tilde{z}(t)+B\tilde{u}(t)+\psi(t)
\end{align}
with
$
\psi(t) = -K_1\nu(t)+\omega(t),
$
and
\begin{align*}
\E \{\psi(t)\psi^T(t)\} &= 
\begin{bmatrix} I & -K_1 \end{bmatrix} 
\begin{bmatrix} \mathcal{W} & \mathcal{U} \\ \mathcal{U}^T& \mathcal{V}\end{bmatrix}
\begin{bmatrix} I & -K_1 \end{bmatrix}^T\\
&= \mathcal{W} -K_1\mathcal {U}^T - \mathcal{U} K_1^T + K_1\mathcal{V} K_1^T.
\end{align*}
Note that $\tilde{z}(t)$ is uncorrelated with $\omega(t), \nu(t)$, since ${z}(t)$ and $\hat{z}(t)$ are uncorrelated with
$\omega(t), \nu(t)$. The gain $K_1$ depends on the dynamics of $\tilde{z}$, which in turn depends on $\tilde{u}$ in (\ref{kalmanz2}), and can be found once we have the dynamics of $\tilde{u}$.

\subsubsection{Cost Decomposition}
Note that 
$$
\begin{bmatrix} x(t) \\ u(t) \end{bmatrix} =
\begin{bmatrix} \hat{x}(t) \\u(t) \end{bmatrix} +
\begin{bmatrix} \tilde{x}(t) \\ 0 \end{bmatrix} 
=
\begin{bmatrix} \hat{z}(t) \\ \hat{u}(t) \end{bmatrix}
+
\begin{bmatrix} \tilde{z}(t) \\ \tilde{u}(t) \end{bmatrix}
+
\begin{bmatrix} \tilde{x}(t) \\ 0 \end{bmatrix} .
$$
The components in the second equality above are uncorrelated, so
\begin{align*}
\mathbf{E} \begin{bmatrix} x(t) \\ u(t)\end{bmatrix}^T \begin{bmatrix}Q&S\\S^T & R \end{bmatrix}  \begin{bmatrix} x(t) \\ u(t) \end{bmatrix}&=\mathbf{E}\begin{bmatrix}  \hat{z}(t) \\  \hat{u}(t) \end{bmatrix}^T \begin{bmatrix}Q&S\\S^T & R \end{bmatrix} \begin{bmatrix}  \hat{z}(t) \\  \hat{u}(t) \end{bmatrix} +\\
\mathbf{E}\begin{bmatrix}  \tilde{z}(t) \\  \tilde{u}(t) \end{bmatrix}^T \begin{bmatrix}Q&S\\S^T & R \end{bmatrix}  \begin{bmatrix}  \tilde{z}(t) \\  \tilde{u}(t) \end{bmatrix} &+  \E~\tilde{x}^T(t) Q\tilde{x}(t) ,
\end{align*}
and hence, we decompose the cost function into three components as $J=\hat{J}_z+\tilde{J}_z+\tilde{J}_x$, with
$$
\hat{J}_z = \lim_{M \rightarrow \infty} \frac{1}{M} \sum_{t=0}^{M-1}\mathbf{E}\begin{bmatrix}  \hat{z}(t) \\  \hat{u}(t) \end{bmatrix}^T \begin{bmatrix}Q&S\\S^T & R \end{bmatrix} \begin{bmatrix}  \hat{z}(t) \\  \hat{u}(t) \end{bmatrix} ,
$$
$$
\tilde{J}_z = \lim_{M \rightarrow \infty} \frac{1}{M} \sum_{t=0}^{M-1}\mathbf{E}\begin{bmatrix}  \tilde{z}(t) \\  \tilde{u}(t) \end{bmatrix}^T \begin{bmatrix}Q&S\\S^T & R \end{bmatrix}  \begin{bmatrix}  \tilde{z}(t) \\  \tilde{u}(t) \end{bmatrix},
$$
and
$$\tilde{J}_x = \lim_{M \rightarrow \infty} \frac{1}{M} \sum_{t=0}^{M-1} \E ~\tilde{x}^T(t) Q\tilde{x}(t).$$ 
$\tilde{J}_x$ is independent of the controller, and hence, nothing can be done about it.
Minimizing $\hat{J}_z$ and $\tilde{J}_z$ is done with respect to independent decision variables, $\hat{u}(t)$ and
$\tilde{u}(t)$, repectively.

\subsection{Finding Controller $\hat{u}$}
Seeking the optimal policy that minimizes the cost $\hat{J}_z$ can be cast as the following optimization problem
\begin{align}
      \min_{\hat{u}}&~\lim_{M \rightarrow \infty} \frac{1}{M}~ \sum_{t=0}^{M-1}  \mathbf{E}~\begin{bmatrix} \hat{z}(t) \\ \hat{u}(t) \end{bmatrix}^T\begin{bmatrix}
                          Q & S \\
                          S^T & R
                        \end{bmatrix}\begin{bmatrix} \hat{z}(t)\\ \hat{u}(t)
    \end{bmatrix}                                                           \nonumber \\
\mbox{subject to}&~~\hat{z}(t+1)=A\hat{z}(t)+B\hat{u}(t)+K_1 (C_1\tilde{z}(t) + \nu(t)).                \label{opt1}
\end{align}
$\hat{u}(t)$ is a function of $y_1(0:t-1)$ by construction. Hence, it can be a function of $\hat{z}(t)$.  This is clearly a centralized state feedback LQ problem. Under stabilizability conditions, the optimal controller is given by $\hat{u}^{*}(t)=-L\hat{z}(t)$, 
where $L$ is obtained from the Riccati equation (\ref{riccati3}) (\cite{c14}).

\subsection{Finding controller $\tilde{u}$}
Note first that $u_1(t)$ depends only on $y_1(0:t-1)$, so  $\hat{u}_1(t) = u_1(t)$. Thus,  $\tilde{u}_1(t) = 0$ and
$$
\tilde{u}(t) = \begin{bmatrix} 0 \\ \tilde{u}_2(t) \end{bmatrix}.
$$

Now the optimal controller synthesis for the cost $\tilde{J}_z$ becomes:
\begin{align*}
      \min_{\tilde{u}_2}&~\lim_{M \rightarrow \infty} \frac{1}{M}~ \sum_{t=0}^{M-1}  \mathbf{E}~\begin{bmatrix} \tilde{z}(t) \\ \tilde{u}_2(t) \end{bmatrix}^T\begin{bmatrix}
                          Q & S_2 \\
                          S_2^T & R_{22}
                        \end{bmatrix}\begin{bmatrix} \tilde{z}(t)\\ \tilde{u}_2(t)
    \end{bmatrix}                                                           \\
\mbox{subject to}&~~\tilde{z}(t+1)=(A-K_1C_1)\tilde{z}(t)+B_2\tilde{u}_2(t)+\psi(t)
\end{align*}
Clearly, this problem is also a centralized state feedback control problem since $\tilde{u}_2$ measures $\tilde{z}$. Therefore, under
stabilizability conditions, the optimal controller is given by
\begin{align}
\tilde{u}_2^{*}(t)&=-L_2\tilde{z}(t),  \label{u2*} 
\end{align}
where $L_2 = R^{-1} B_2^T \Pi_2$, and $\Pi_2$ is given by the Riccati equation (\ref{riccati3}) (\cite{c14}).

\subsection{Finding the Estimator Gain $K_1$}
First, take the optimal controller (\ref{u2*}) and plug it in the expression of the Kalman filter
error equation (\ref{kalmanz2}):
\begin{align}\label{kalmanclosed}
    \tilde{z}(t+1)=(A-K_1C_1-B_2L_2)\tilde{z}(t)+\psi(t).
\end{align}
Thus, the error dynamics of $\tilde{z}$ give rise to the Riccati equation (\ref{riccati6}) (\cite{c14}),
and $K_1$ is the solution obtained from  (\ref{riccati6}) and (\ref{K2}).

\section{Discussion}
We can see that the controller is decomposed into two parts,
$$
\begin{bmatrix} {u}_1(t) \\ u_2(t) \end{bmatrix} = \begin{bmatrix} \hat{u}_1(t) \\ \hat{u}_2(t) \end{bmatrix} + 
\begin{bmatrix} 0 \\ \tilde{u}_2(t) \end{bmatrix},
$$
where 
$$
\begin{bmatrix} \hat{u}_1(t) \\ \hat{u}_2(t) \end{bmatrix} 
$$
is a state feedback law based on $\hat{z}(t)$. $\hat{z}(t)$ is an estimate of the \emph{centralized} estimate of the state $x(t)$, based on the output history $y_1(0:t-1)$ which is available to both controllers. The additional
information controller 2 has access to, $y_2(0:t-1)$, is used to to minimize the part of the cost ($\tilde{J}_z$) that can not be observed with respect to $y_1(0:t-1)$. This is given by $ \tilde{u}_2(t)$.

The dependence of the filter (\ref{kalmanz2}) on part of the controller, $\tilde{u}_2$, is essentially what makes the problem different from
centralized control, where the estimator becomes independent of the controller. This is also why we get a coupling between
the estimator gain $K_1$ and the controller gain $L_2$ in the Riccati equations  (\ref{riccati6})-(\ref{riccati5}).

We note that the Riccati equations (\ref{riccati4})-(\ref{riccati3}) can be solved independently.
The Riccati equations (\ref{riccati6})-(\ref{riccati5}) are coupled and depend on the output from (\ref{riccati4})-(\ref{riccati3}). 
The coupling looks complicated, but using the fact that the equalities in Riccati equations can be written as inequalities (\cite{c15}) , 
and by Schur complementing, we can obtain a set of linear matrix inequalities that can be efficiently solved. 

Also, in principle, the solution can be easily extended to the case of $N$ interconnected
systems by induction. We simply take the $N$th system as the "second system" considerd in this paper, and the rest of 
systems $1, ..., N-1$ to be the "first" system. Then we proceed iteratively by taking system $N-1$ as the "second" system and 
systems $1, ..., N-2$ as the first, and so on. However, the coupled Riccati equations might give rise to a complicated numerical solution.

\end{document}